\documentclass[11pt]{amsart}
\usepackage[margin=1.05in]{geometry}
\usepackage{amsmath,amssymb,amsthm,mathtools}
\usepackage{enumitem}
\usepackage[normalem]{ulem}
\usepackage{mathrsfs}
\usepackage[dvipsnames]{xcolor}
\usepackage[colorlinks=true,linkcolor=Plum,citecolor=Maroon,urlcolor=SpringGreen]{hyperref}
\newtheorem{theorem}{Theorem}[section]
\newtheorem{proposition}[theorem]{Proposition}
\newtheorem{lemma}[theorem]{Lemma}
\newtheorem{corollary}[theorem]{Corollary}
\theoremstyle{definition}
\newtheorem{definition}[theorem]{Definition}
\newtheorem{remark}[theorem]{Remark}
\newtheorem{example}[theorem]{Example}
\newtheorem{question}[theorem]{Question}
\numberwithin{equation}{section}
\newcommand{\PP}{\mathbb{P}}
\newcommand{\CC}{\mathbb{C}}
\newcommand{\ZZ}{\mathbb{Z}}
\newcommand{\QQ}{\mathbb{Q}}
\newcommand{\cO}{\mathcal{O}}
\newcommand{\cF}{\mathcal{F}}
\newcommand{\cI}{\mathcal{I}}
\newcommand{\cJ}{\mathcal{J}}
\newcommand{\cA}{\mathcal{A}}
\newcommand{\cD}{\mathcal{D}}
\newcommand{\sD}{\mathscr{D}}
\newcommand{\Pfrak}{\mathfrak{P}}
\newcommand{\Pcomb}{\mathfrak{P}_{\rm comb}}
\newcommand{\pt}{\mathrm{pt}}
\newcommand{\tb}{\widetilde{b}}
\newcommand{\tR}{\widetilde{R}}
\newcommand{\talpha}{\widetilde{\alpha}}
\newcommand{\gadd}[1]{#1}
\newcommand{\gdel}[1]{}
\newcommand{\gblock}{}
\DeclareMathOperator{\AR}{AR}
\DeclareMathOperator{\Sing}{Sing}
\DeclareMathOperator{\mdr}{mdr}
\DeclareMathOperator{\length}{length}
\DeclareMathOperator{\lct}{lct}
\DeclareMathOperator{\lcm}{lcm}
\DeclareMathOperator{\mult}{mult}

\title{Chern classes, Poincar\'e polynomials, and Bernstein--Sato polynomials of plane curves}
\date{\today}
\author{Piotr Pokora}
\begin{document}

\begin{abstract}
We show that the Poincar\'e polynomial of a reduced plane curve $C$ is the Chern polynomial
of the dual of the bundle of Jacobian syzygies of $C$. Using the Zeuthen--Segre formula,
we prove that its value at $t=-1$ equals the Euler number of the complement of $C$ up to
the correction term $\mu(C)-\tau(C)$, so the two numbers coincide if and only if all
singularities of $C$ are quasi-homogeneous. As applications, we show that the defect in a
combinatorial inequality for free $d$-arrangements is the square of the difference of the
exponents, extend the addition technique to arbitrary reduced plane curves, and locate the
roots of local Bernstein--Sato polynomials at ordinary singular points. 
\end{abstract}
\maketitle
\section{Introduction}

In the present paper we continue the studies on Poincar\'e polynomials of reduced plane curves initiated in \cite{Di25, Po25, Po24}.
Let us briefly recall the main object of our considerations.

\begin{definition}
Let $C\colon f=0$ be a reduced curve in $\PP^2_\CC$ of degree $d$.
The Poincar\'e polynomial of $C$ is defined as
\[
\Pfrak(C;t)=1+(d-1)t+\bigl((d-1)^2-\tau(C)\bigr)t^2,
\]
where $\tau(C)$ denotes the total Tjurina number of $C$.
\end{definition}

It is known that if $C$ is free with exponents $(d_1,d_2)$, then $\Pfrak(C;t)=(1+d_1t)(1+d_2t)$, see \cite{Di25, Po24}.
Moreover, it was observed in \cite[Remark~3.5]{Po25} that for conic-line arrangements with ordinary \gdel{quasi-homogeneous} singularities the value of the combinatorial Poincar\'e polynomial at $t=-1$ is equal to the Euler number of the complement $M(C):=\PP^2_\CC\setminus C$.
It is natural to wonder whether there is a natural and conceptual explanation of these phenomena.
The main aim of the present paper is to provide such an explanation by using Chern classes and exact sequences, following the approach to the Zeuthen--Segre formula presented by Catanese in \cite[Section~1.3.3]{BCP}.

Let us denote by $E_C$ the sheaf on $\PP^2_\CC$ associated with the graded module of Jacobian syzygies of $f$. 

Our first main result gives a clear relation between the Euler characteristic of $M(C)$ and the Poincar\'e polynomial of $C$ in the language of the Tjurina and Milnor numbers.

\begin{theorem}\label{thm:A}
Let $C\subset\PP^2_\CC$ be a reduced curve. Then $\Pfrak(C;t)=c_t(E_C^\vee)$, and
\[
e(M(C))=\Pfrak(C;-1)-\mu(C)+\tau(C),
\]
where $\mu(C)$ denotes the total Milnor number of $C$.
In particular, $e(M(C))=\Pfrak(C;-1)$ if and only if $C$ admits only quasi-homogeneous singularities.
\end{theorem}

In other words, the Poincar\'e polynomial is a Chern polynomial, and its value at $t=-1$ computes the Euler number of $M(C)$ up to the correction term $\mu(C)-\tau(C)$, which vanishes exactly in the quasi-homogeneous case.
Since free curves are exactly those curves for which $E_C$ splits, we immediately obtain the following corollary, which is also a special case of \cite[Theorem 2.5 iv]{DS17}.

\begin{corollary}\label{cor:A}
Let $C\subset\PP^2_\CC$ be a free curve with exponents $(d_1,d_2)$ admitting only quasi-homogeneous singularities. Then
\[
e(M(C))=(d_1-1)(d_2-1).
\]
\end{corollary}

Next, recall that an arrangement $\cA=\{C_1,\dots,C_k\}\subset\PP^2_\CC$ is called a $d$-arrangement if all curves $C_i$ are smooth of the same degree $d\ge1$ and all intersection points of $\cA$ are ordinary singularities.
By \cite[Theorem~1.8]{Po25}, if $\cA$ is a free $d$-arrangement with $k\ge2$ and $d\ge2$ having only quasi-homogeneous singularities, then
\[
\sum_{r\ge2}(r^2-5r+4)n_r\ge3+3dk(d-2),
\]
where $n_r$ denotes the number of $r$-fold points.
Our second result shows that the difference of both sides is always a perfect square, and this result clarifies the meaning of our numerical criterion for free $d$-arrangements \cite[Theorem~1.8]{Po25}.

\begin{theorem}\label{thm:B}
Let $\cA$ be a $d$-arrangement with $k\ge2$ and $d\ge2$.
Assume that all singularities of $\cA$ are quasi-homogeneous and that $\cA$ is free with exponents $(d_1,d_2)$. Then
\[
\sum_{r\ge2}(r^2-5r+4)n_r-3-3dk(d-2)=(d_2-d_1)^2 .
\]
\end{theorem}

Then we revisit the addition technique that was introduced in \cite{Po24}, and to approach this problem we will use the following idea introduced by Dimca in \cite{Di24}. For a reduced curve $C$ we put $\varepsilon(C) :=\mu(C)-\tau(C)\ge0$.
Then for two reduced curves $C_1,C_2$ without common components we put 
\[
\varepsilon(C_1,C_2) := \varepsilon(C_1\cup C_2)-\varepsilon(C_1)-\varepsilon(C_2).
\]
These two notations allow us to formulate an addition technique for arbitrary plane curves in the complex projective plane that have no common components.
\begin{theorem}\label{thm:C}
Let $C_1,C_2\subset\PP^2_\CC$ be reduced curves without common components and let $r=|C_1\cap C_2|$. Then
\[
\Pfrak(C_1\cup C_2;t)=\Pfrak(C_1;t)+\Pfrak(C_2;t)+t-1+\bigl(r-1+\varepsilon(C_1,C_2)\bigr)t^2 .
\]
In particular, the addition formula \cite[Theorem~1.2]{Po24} holds provided that $C_1\cup C_2$ is quasi-homogeneous at all points of $C_1\cap C_2$.
\end{theorem}

In the final part of the paper, we reformulate our results on free plane curves in the language of Bernstein–Sato polynomials. This direction is a natural one: the log canonical threshold, which plays an important role in the homological study of plane curves, is determined by the largest root of the Bernstein–Sato polynomial.

Let us recall that for a point $p$ of a reduced curve $C$ we denote by $b_{C,p}(s)$ the local Bernstein--Sato polynomial of a local equation of $C$ at $p$, and by $R_p(C)$ the set of roots of $b_{C,p}(-s)$.
By Kashiwara's theorem \cite{Ka}, $R_p(C)$ is a finite subset of $\QQ_{>0}$.
It turns out that quasi-homogeneous and ordinary singularities, which play the main role in Theorems~\ref{thm:A}--\ref{thm:C}, are also singularities for which the set $R_p(C)$ can be effectively controlled.
Recall that $p\in C$ is an ordinary $r$-fold point if $C$ has at $p$ exactly $r$ smooth branches which are pairwise transverse.
Our last main result is the following characterization. We believe that this result might be well-known to the experts, but here we present it in the version adapted to our considerations.

\begin{theorem}\label{thm:D}
Let $C\subset\PP^2_\CC$ be a reduced curve and let $p\in C$ be an ordinary $r$-fold point with $r\ge2$. Then
\[
\Bigl\{\frac jr\colon 2\le j\le \gadd{\min\{r+1,2r-2\}}\Bigr\}\subseteq R_p(C)\subseteq\Bigl\{\frac jr\colon 2\le j\le 2r-2\Bigr\}.
\]
Moreover, if $p$ is quasi-homogeneous, then
\[
b_{C,p}(s)=(s+1)\prod_{j=2}^{2r-2}\Bigl(s+\frac jr\Bigr),
\]
so the second inclusion is an equality.
\end{theorem}

Combining Theorems~\ref{thm:B} and~\ref{thm:D}, we obtain the following result, where $b_\cA(s)$ denotes the Bernstein--Sato polynomial of $\cA$, i.e., the least common multiple of the local Bernstein--Sato polynomials $b_{\cA,p}(s)$.

\begin{corollary}\label{cor:D}
Let $\cA$ be a free $d$-arrangement with $k\ge2$ and $d\ge2$ having only quasi-homogeneous singularities.
Then the maximal multiplicity $m$ of a singular point of $\cA$ satisfies $5\le m\le k$, \gadd{the log canonical threshold of $\cA$ is $2/m\le 2/5$,} and $b_\cA(s)$ has at least $2m-3\ge7$ pairwise distinct roots.
\end{corollary}

We also show that in the quasi-homogeneous case the roots of the reduced local Bernstein--Sato polynomials are symmetric with respect to $-1$, and that the mechanism behind Theorem~\ref{thm:C} leads to an additivity property of log canonical thresholds that seems to be not explicitly stayed in the literature. Finally, we formulate a question which, if answered positively, would give a characterization of the equality $e(M(C))=\Pfrak(C;-1)$ purely in terms of Bernstein--Sato polynomials.

Let us present a short outline of our paper.

In Section~\ref{sec:prel} we recall basic facts on free plane curves, logarithmic vector fields, and Chern classes.
In Section~\ref{sec:A} we prove Theorem~\ref{thm:A} and Corollary~\ref{cor:A}.
In Section~\ref{sec:comb} we explain how the combinatorial Poincar\'e polynomials of conic-line arrangements from \cite{Po25} fit into this picture.
In Section~\ref{sec:B} we prove Theorem~\ref{thm:B}.
In Section~\ref{sec:C} we prove Theorem~\ref{thm:C} and we discuss the addition technique in the language of logarithmic vector fields.
Finally, in Section~\ref{sec:BS} we prove Theorem~\ref{thm:D} and Corollary~\ref{cor:D}, and we discuss further relations between our results and Bernstein--Sato polynomials.

We work exclusively over the complex numbers.

\section{Preliminaries}\label{sec:prel}

We gather all necessary tools that we are going to use in the paper, and in our presentation we are going to follow the notation introduced in \cite{Di17}.

Let $S:=\CC[x,y,z]$ be the coordinate ring of $\PP^2_\CC$ and let $C\colon f=0$ be a reduced curve of degree $d$.
We denote by $J_f=\langle\partial_xf,\partial_yf,\partial_zf\rangle$ the Jacobian ideal of $f$. Let $p$ be a singular point of $C$ and let $h$ be a local equation of $C$ at $p$ in local analytic coordinates $(u,v)$ centred at $p$.
The Milnor and the Tjurina number of $C$ at $p$ are defined as
\[
\mu_p=\dim_\CC\cO_p/\langle\partial_uh,\partial_vh\rangle\qquad\text{and}\qquad\tau_p=\dim_\CC\cO_p/\langle h,\partial_uh,\partial_vh\rangle,
\]
and we set 
$$\mu(C)=\sum_p\mu_p \text{ and }\tau(C)=\sum_p\tau_p,$$ where the sums run over all singular points of $C$.
The algebras $\cO_p/\langle\partial_uh,\partial_vh\rangle$ and $\cO_p/\langle h,\partial_uh,\partial_vh\rangle$ are called the Milnor and the Tjurina algebra of $C$ at $p$, respectively.
By the constuction we know that $\tau_p\le\mu_p$.
We say that $p$ is a quasi-homogeneous singular point if there are local analytic coordinates centred at $p$ in which $C$ is defined by a weighted homogeneous polynomial.
By a result of K.~Saito \cite{Sa71}, $p$ is quasi-homogeneous if and only if $\tau_p=\mu_p$.

A point $p\in C$ of multiplicity $r\ge2$ is called an $r$-fold point.
We say that an $r$-fold point $p$ is ordinary if $C$ has at $p$ exactly $r$ smooth branches which are pairwise transverse.
Equivalently, the tangent cone of $C$ at $p$ consists of $r$ pairwise distinct lines.
By Milnor's formula, one has $\mu_p=2\delta_p-\gadd{r_p}+1$\gadd{, where $r_p$ is the number of branches of $C$ at $p$}, and for an ordinary $r$-fold point \gadd{one has $r_p=r$ and $\delta_p=\binom r2$, which} gives us $\mu_p=(r-1)^2$.
Ordinary $r$-fold points with $r\le4$ are always quasi-homogeneous, being of type $A_1$, $D_4$ or $X_9$, but for $r\ge5$ this is no longer the case, see Example~\ref{ex:nonqh}.

\begin{definition}
Consider the graded $S$-module of Jacobian syzygies of $f$, namely
\[
\AR(f)=\{(a,b,c)\in S^3\colon a\partial_xf+b\partial_yf+c\partial_zf=0\},
\]
and let $\mdr(f):=\min\{r\ge0\colon\AR(f)_r\neq0\}$ be the minimal degree of a Jacobian relation for $f$.

A reduced curve $C\colon f=0$ of degree $d$ is \textbf{free} if $\AR(f)$ is a free $S$-module, or equivalently, if $J_f$ is saturated with respect to $\langle x,y,z\rangle$. In this case $\AR(f)\cong S(-d_1)\oplus S(-d_2)$ with $d_1=\mdr(f)$ and $d_1+d_2=d-1$, and the pair of non-negative integers $(d_1,d_2)$ is called the exponents of $C$.
\end{definition}

The following result is crucial for our considerations.

\begin{theorem}[du Plessis--Wall, \cite{dPW}]\label{thm:dPW}
Let $C\colon f=0$ be a reduced curve in $\PP^2_\CC$ of degree $d$ and let $d_1=\mdr(f)$.
Then $C$ is free if and only if $(d-1)^2-d_1(d-d_1-1)=\tau(C)$.
\end{theorem}

We denote by $E_C$ the sheaf on $\PP^2_\CC$ associated with the graded $S$-module $\AR(f)$.
Since $\AR(f)$ is the kernel of a morphism of graded free $S$-modules, we have $\AR(f)=\bigoplus_{k\in\ZZ}H^0(\PP^2_\CC,E_C(k))$, and hence $C$ is free if and only if $E_C$ splits as a direct sum of line bundles.

Let $T\langle C\rangle$ denote the sheaf of logarithmic vector fields along $C$,
that is, the sheaf of vector fields $\theta$ on $\PP^2_\CC$ such that
$\theta(h) \in (h)$ for every local equation $h$ of $C$.
It is well-known that $T\langle C\rangle \cong E_C(1)$, see for instance \cite{DS}. In particular,
if $C$ is free with exponents $(d_1, d_2)$, then
\[
  T\langle C\rangle \cong \cO(1-d_1) \oplus \cO(1-d_2).
\]

For a coherent sheaf $\mathcal E$ on $\PP^2_\CC$ we write its Chern polynomial $c_t(\mathcal E)=1+c_1(\mathcal E)t+c_2(\mathcal E)t^2$, where we identify $H^2(\PP^2_\CC,\ZZ)$ and $H^4(\PP^2_\CC,\ZZ)$ with $\ZZ$ via the hyperplane class $h$ and the point class.
Let us recall the following facts, which are the main ingredients of Catanese's proof of the Zeuthen--Segre formula \cite[Theorem~1.33]{BCP}.
\begin{enumerate}[label=\textup{(C\arabic*)},leftmargin=*]
\item\label{C1} The total Chern class is multiplicative on short exact sequences of coherent sheaves.
\item\label{C2} If $Q$ is a coherent sheaf of finite length $\ell$ on a smooth
projective surface, then $Q\otimes L\cong Q$ for every line bundle $L$. In particular,
$c(Q\otimes L)=c(Q)$. Moreover, $c(Q)=1-\ell\,[\pt]$ in cohomology
(or in the Chow ring, if the surface is rational, as in our case).
\item\label{C3} If $\varphi\colon X\to B$ is a fibration of a smooth projective surface onto a smooth curve with general fibre $F$, and $\Sigma\subset B$ is the set of critical values, then
\[
e(X)=e(B)e(F)+\sum_{y\in\Sigma}\bigl(e(F_y)-e(F)\bigr),
\]
see for instance \cite[Chapter~III, Section~11]{BHPV}.
\end{enumerate}

The following special case of \cite[Theorem~1.33]{BCP} is all that we need here.

\begin{theorem}[Zeuthen--Segre]\label{thm:ZS}
Let $\varphi\colon X\to B$ be a fibration of a smooth projective surface onto a smooth curve such that all fibres are reduced, and let $F$ be a general fibre.
Then $e(X)=e(B)e(F)+\deg\cF$, where $\cF$ is the structure sheaf of the critical scheme of $\varphi$.
Locally, $\cF=\cO_X/\langle\partial s/\partial x_1,\partial s/\partial x_2\rangle$, where $s$ is the pull-back of a local parameter of $B$.
In particular, the length of $\cF$ at a singular point $p$ of a fibre $F_y$ is equal to the Milnor number of $F_y$ at $p$.
\end{theorem}

\begin{proof}
Since all fibres are reduced, the exact sequence $(*)$ from \cite[p.~28]{BCP} takes the form
\[
0\to\varphi^*\omega_B\to\Omega^1_X\to\omega_{X/B}\to\cF\otimes\omega_{X/B}\to0 .
\]
By \ref{C1} and \ref{C2}, we have
\[
c(\Omega^1_X)=(1+\varphi^*K_B)(1+K_X-\varphi^*K_B)(1+\deg\cF\,[\pt]).
\]
Since $(\varphi^*K_B)^2=0$, $\varphi^*K_B\equiv(2b-2)F$ and $K_X\cdot F=2g(F)-2$, where $b$ is the genus of $B$, we get
\[
e(X)=c_2(\Omega^1_X)=(2b-2)(2g(F)-2)+\deg\cF=e(B)e(F)+\deg\cF,
\]
which completes the proof.
\end{proof}

\section{Poincar\'e polynomials via Chern classes}\label{sec:A}

Let $C\colon f=0$ be a reduced curve of degree $d$ and let $E_C$ be the sheaf associated with $\AR(f)$.
We start with the following well-known fact.

\begin{lemma}\label{lem:syz}
Let $Z\subset\PP^2_\CC$ be the zero-dimensional scheme defined by $J_f$. Then there is an exact sequence
\begin{equation}\label{eq:syz}
0\to E_C\to\cO^3_{\PP^2_\CC}\to\cI_Z(d-1)\to0,
\end{equation}
so $E_C$ is a rank two vector bundle, and $\deg Z=\tau(C)$.
\end{lemma}
Let us now express our Poincar\'e polynomial as the Chern polynomial of $E_{C}^{\vee}$.
\begin{proposition}\label{prop:chern}
One has $c(E_C)=1-(d-1)h+\bigl((d-1)^2-\tau(C)\bigr)h^2$, and hence $\Pfrak(C;t)=c_t(E_C^\vee)$.
\end{proposition}

\begin{proof}
By \ref{C1} and \ref{C2} applied to the exact sequence 
$$0\to\cI_Z(d-1)\to\cO(d-1)\to\cO_Z\to 0,$$
we get 
$$c(\cI_Z(d-1))=1+(d-1)h+\tau(C)h^2 .$$
Now \ref{C1} applied to \eqref{eq:syz} gives
\[
c(E_C)=c(\cI_Z(d-1))^{-1}=\bigl(1+(d-1)h+\tau(C)h^2\bigr)^{-1}=1-(d-1)h+\bigl((d-1)^2-\tau(C)\bigr)h^2 .
\]
Since $c_1(E_C^\vee)=-c_1(E_C)$ and $c_2(E_C^\vee)=c_2(E_C)$, we obtain $\Pfrak(C;t)=c_t(E_C^\vee)$, which completes the proof.
\end{proof}

\begin{remark}\label{rem:split}
If $C$ is free with exponents $(d_1,d_2)$, then $E_C\cong\cO(-d_1)\oplus\cO(-d_2)$, so the splitting $\Pfrak(C;t)=(1+d_1t)(1+d_2t)$ is immediate.
Comparing the second Chern classes, we recover the equality $\tau(C)=(d-1)^2-d_1d_2$ from Theorem~\ref{thm:dPW}.
\end{remark}

Recall that $T\langle C\rangle\cong E_C(1)$.
Since $c_2(V(1))=c_2(V)+c_1(V)+1$ for a rank two vector bundle $V$, Proposition~\ref{prop:chern} gives
\begin{equation}\label{eq:twist}
\Pfrak(C;-1)=1+c_1(E_C)+c_2(E_C)=c_2(T\langle C\rangle)=d^2-3d+3-\tau(C).
\end{equation}
In other words, evaluating the Poincar\'e polynomial at $t=-1$ corresponds to the twist by $\cO(1)$.

Our next step is to compute the Euler number of a reduced plane curve by means of
Theorem~\ref{thm:ZS}. This is a well-known result. By \cite[Section~9.2, Lemma~1 and the
proof of Theorem~5]{BK}, one has $e(C)=e(C_t)-\sum_{p}\bigl(e(F_p)-1\bigr)$, where $C_t$ is
a smooth plane curve of degree $d$ and $F_p$ is the Milnor fibre of $C$ at $p$, and
$e(F_p)=1-\mu_p$ by \cite[Section~8.5]{BK}. Here we present a different proof, based on the
Zeuthen--Segre theorem.

\begin{proposition}\label{prop:eC}
Let $C\subset\PP^2_\CC$ be a reduced curve of degree $d$. Then $e(C)=3d-d^2+\mu(C)$.
\end{proposition}

\begin{proof}
We may assume that $d\ge2$.
Let $\cD\subset\gadd{\PP^n}:=\PP(H^0(\cO(d)))$ be the discriminant hypersurface, and let $\cD^\circ\subset\cD$ be the dense open subset of curves with exactly one singular point which is an ordinary node.

\emph{The pencil.}
Take a general $g\in H^0(\cO(d))$.
Then $G\colon g=0$ is smooth and meets $C$ transversally at $d^2$ smooth points of $C$.
Since $\dim(\cD\setminus\cD^\circ)\le \gadd{n}-2$, the cone over $\cD\setminus\cD^\circ$ with vertex $[f]$ has dimension at most $\gadd{n}-1$, so the line through $[f]$ and $[g]$ meets $\cD\setminus\cD^\circ$ at most at $[f]$.
Hence every member of the pencil $\lambda_0f+\lambda_1g=0$ different from $C$ is either smooth or has exactly one singular point, which is an ordinary node.
Moreover, every member of the pencil is smooth at the $d^2$ base points, since $C$ and $G$ are smooth and transverse there.
Let $N$ denote the number of nodal members.

\emph{The fibration.}
Blowing up the base points, we obtain a fibration $\varphi\colon X\to\PP^1$ whose fibres are reduced and isomorphic to the members of the pencil.
We have $e(X)=3+d^2$ and $e(F)=2-(d-1)(d-2)=3d-d^2$ for a general fibre $F$.
Near a singular point $p$ of $C$ the pull-back of a local parameter of $\PP^1$ is $s=f/g=u\cdot h$, where $u$ is a unit and $h$ is a local equation of $C$.
Since the Milnor number of a reduced plane curve germ depends only on its zero set \cite[Section~10]{Mi}, the length of $\cF$ at $p$ is equal to $\mu_p$.

\emph{Comparison.}
Theorem~\ref{thm:ZS} gives
\[
3+d^2=2(3d-d^2)+\mu(C)+N .
\]
On the other hand, \ref{C3} together with the equality $e(C')=e(F)+1$ for a nodal member $C'$, which follows by normalization, gives
\[
3+d^2=2(3d-d^2)+\bigl(e(C)-e(F)\bigr)+N .
\]
Comparing both identities, we get $e(C)=e(F)+\mu(C)=3d-d^2+\mu(C)$, which completes the proof.
\end{proof}

As a by-product, we obtain $N=3(d-1)^2-\mu(C)$, and recall that for smooth $C$ this is the classical number of singular members of a general pencil of plane curves of degree $d$.
Proposition~\ref{prop:eC} can also be obtained topologically from Milnor's formula $\mu_p=2\delta_p-r_p+1$, where $r_p$ is the number of branches of $C$ at $p$.

Now we are ready to prove our first main result.

\begin{proof}[Proof of Theorem~\ref{thm:A}]
The first part is Proposition~\ref{prop:chern}.
By Proposition~\ref{prop:eC} and \eqref{eq:twist}, we have
\[
e(M(C))=e(\PP^2_\CC)-e(C)=d^2-3d+3-\mu(C)=\Pfrak(C;-1)-\mu(C)+\tau(C).
\]
Since $\tau_p\le\mu_p$ for every singular point $p$, with equality if and only if $p$ is quasi-homogeneous \cite{Sa71}, we have $\mu(C)=\tau(C)$ if and only if $C$ admits only quasi-homogeneous singularities, which completes the proof.
\end{proof}

\begin{proof}[Proof of Corollary~\ref{cor:A}]
By Theorem~\ref{thm:A} and Remark~\ref{rem:split}, we have $e(M(C))=\Pfrak(C;-1)=(1-d_1)(1-d_2)=(d_1-1)(d_2-1)$.
\end{proof}

\begin{remark}
Note that $d_1=\mdr(f)=0$ if and only if $C$ is a pencil of $d$ lines, and then $e(M(C))=2-d$.
In all remaining cases Corollary~\ref{cor:A} gives $e(M(C))\ge0$, with equality if and only if $d_1=1$.
\end{remark}

\begin{remark}
Theorem~\ref{thm:A} can be compared with the results of Aluffi \cite{Al} and Liao \cite{Li} on the Chern--Schwartz--MacPherson class of the complement and the Chern class of the sheaf of logarithmic vector fields.
\end{remark}

\begin{example}\label{ex:freeA}
Consider the arrangement $\cA$ of $9$ lines and one conic that comes from \cite[Example~3.3]{Po25}. We know that this arrangement is free with exponents $(4,6)$ and satisfies $\tau(\cA)=\mu(\cA)=76$.
Corollary~\ref{cor:A} gives $e(M(\cA))=3\cdot5=15$.
Since $\cA$ has only ordinary singularities with $(n_2,n_3,n_4)=(6,4,6)$, the direct count from Proposition~\ref{prop:comb} below gives $e(M(\cA))=3-2\cdot10+(6+8+18)=15$, as expected.
Similarly, for the free arrangement of three conics from \cite[Example~3.9]{Po25} with exponents $(2,3)$ \gadd{and only quasi-homogeneous singularities} we get $e(M(\cA))=1\cdot2=2$.
\end{example}

\begin{example}\label{ex:nonqh}
Let us consider the sextic curve
\[
C\colon f=x^5z+y^5z+x^3y^3=0 .
\]
We claim that $p=[0:0:1]$ is the only singular point of $C$.
Indeed, the line $z=0$ meets $C$ at $[1:0:0]$ and $[0:1:0]$, which are smooth points of $C$ since $\partial_zf=x^5+y^5$ does not vanish there.
In the affine chart $z=1$ we put $g=x^5+y^5+x^3y^3$, and we have $x\partial_xg+y\partial_yg=5g+x^3y^3$.
Hence at a singular point we have $xy=0$, and then $\partial_xg=\partial_yg=0$ gives $x=y=0$.
The point $p$ is an ordinary $5$-fold point, since the tangent cone $x^5+y^5=0$ consists of five pairwise distinct lines, so $\mu(C)=16$.
Using \texttt{SINGULAR} script from \textsl{Appendix}, we can check that $\tau(C)=\tau_p=15$, so $p$ is not quasi-homogeneous.
Hence
\[
\Pfrak(C;t)=1+5t+10t^2,\qquad \Pfrak(C;-1)=6,\qquad e(M(C))=36-18+3-16=5.
\]
Since $(d-1)^2-d_1(d-1-d_1)\neq15$ for every integer $d_1$, the curve $C$ is not free by Theorem~\ref{thm:dPW}.
Moreover, $C$ is irreducible.
Indeed, if $C=C'\cup C''$ with $\deg C'=a$ and $\deg C''=b$, then $C'\cap C''=\{p\}$, since every point of $C'\cap C''$ is a singular point of $C$, and by B\'ezout's theorem $ab=i_p(C',C'')=s(5-s)$, where $s$ is the number of branches of $C'$ at $p$, which is impossible since $ab\in\{5,8,9\}$ and $s(5-s)\in\{0,4,6\}$.
Hence $C$ is a rational curve, since its geometric genus is equal to $10-\binom52=0$, and $e(C)=-2$ by Proposition~\ref{prop:eC}.
\end{example}

\section{Combinatorial Poincar\'e polynomials}\label{sec:comb}

In this section we explain how the combinatorial Poincar\'e polynomials of conic-line arrangements introduced in \cite{Po25} fit into the picture from Section~\ref{sec:A}.
Let $\cA$ be an arrangement of $d$ lines and $k$ smooth conics with only ordinary singularities.
We put $f_i=\sum_{r\ge2}r^in_r$.
Recall that the combinatorial Poincar\'e polynomial of $\cA$ is defined in \cite{Po25} as
\[
\Pcomb(\cA;t)=1+(2k+d-1)t+\Bigl(\sum_{r\ge2}(r-1)n_r-d+1\Bigr)t^2 .
\]

\begin{proposition}\label{prop:comb}
Let $\cA$ be an arrangement of $d$ lines and $k$ smooth conics with only ordinary singularities. Then
\[
e(M(\cA))=3-2(d+k)+\sum_{r\ge2}(r-1)n_r=\Pcomb(\cA;-1).
\]
\end{proposition}

\begin{proof}
The first equality is a well-known result, see for instance \cite[Section 1.2]{ST}, and the second equality follows directly from the definition of $\Pcomb(\cA;t)$.
\end{proof}

\begin{corollary}\label{cor:comb}
In the setting of Proposition~\ref{prop:comb}, one has
\[
\Pcomb(\cA;t)=\Pfrak(\cA;t)-\bigl(\mu(\cA)-\tau(\cA)\bigr)t^2 .
\]
\end{corollary}

\begin{proof}
Since every $r$-fold point has Milnor number $(r-1)^2$, we have $\mu(\cA)=f_2-2f_1+f_0$.
The na\"ive combinatorial count $4\binom k2+2kd+\binom d2=\sum_r\binom r2n_r$ gives $f_2-f_1=(2k+d-1)^2+d-1$.
Hence $\sum_r(r-1)n_r-d+1=(2k+d-1)^2-\mu(\cA)$, while the coefficient of $t^2$ in $\Pfrak(\cA;t)$ is equal to $(2k+d-1)^2-\tau(\cA)$\gdel{l}, which completes the proof.
\end{proof}

\begin{remark}
Proposition~\ref{prop:comb} and Corollary~\ref{cor:comb} explain the different hypotheses in \cite{Po25}.
The polynomial $\Pcomb$ is of topological nature, since its value at $t=-1$ is the second Chern class of the logarithmic cotangent bundle on a log resolution.
The polynomial $\Pfrak$ is of algebraic nature, since its value at $t=-1$ is the second Chern class of $T\langle\cA\rangle$ on $\PP^2_\CC$, see \eqref{eq:twist}.
Both polynomials coincide exactly in the quasi-homogeneous case.
This is why \cite[Theorem~1.5]{Po25} requires quasi-homogeneity, whereas \cite[Remark~3.5]{Po25} does not.
Since ordinary points of multiplicity $r\le4$ are always quasi-homogeneous, the quasi-homogeneity assumption in \cite[Theorem~1.5]{Po25} concerns only points of multiplicity $r\ge5$.
\end{remark}

\section{Proof of Theorem~\ref{thm:B}}\label{sec:B}

For a rank two vector bundle $E$ on $\PP^2_\CC$ let us define by $\Delta(E)=c_1(E)^2-4c_2(E)$ its discriminant.
By Proposition~\ref{prop:chern}, $\Delta(E_C)$ is the discriminant of the quadratic polynomial $\Pfrak(C;t)$, namely
\[
\Delta(E_C)=4\tau(C)-3(d-1)^2 .
\]
If $C$ is free with exponents $(d_1,d_2)$, then $\Delta(E_C)=(d_1+d_2)^2-4d_1d_2=(d_2-d_1)^2$.

\begin{proof}[Proof of Theorem~\ref{thm:B}]
Keeping the notation of Section~\ref{sec:comb}, for a $d$-arrangement $\cA$ we have $\deg\cA=dk$ and $\tau(\cA)=\mu(\cA)=f_2-2f_1+f_0$.
Since any two curves of $\cA$ meet transversally at $d^2$ points, the combinatorial count $d^2\binom k2=\sum_r\binom r2n_r$ gives $f_2-f_1=d^2(k^2-k)$.
Then
\[
\Delta(E_\cA)=4f_2-8f_1+4f_0-3\bigl(f_2-f_1+d^2k\bigr)+6dk-3=\sum_{r\ge2}(r^2-5r+4)n_r-3-3dk(d-2).
\]
Since $\cA$ is free, we have $\Delta(E_\cA)=(d_2-d_1)^2$, which completes the proof.
\end{proof}

\begin{remark}\label{rem:disc}
The discriminant $\triangle_{d_1}$ in the proof of \cite[Theorem~1.8]{Po25} is precisely $\Delta(E_\cA)$.
Moreover, $d_2-d_1\equiv dk-1\pmod 2$, since $d_1+d_2=dk-1$.

In the same way, for an arrangement $\cA$ of $d$ lines and $k$ smooth conics with ordinary quasi-homogeneous singularities one gets
\[
\Delta(E_\cA)=\sum_{r\ge2}(r^2-5r+4)n_r+3(d-1),
\]
which must be a perfect square if $\cA$ is free.
For $d=0$ this sharpens \cite[Remark~3.11]{Po25}: if a conic arrangement with ordinary quasi-homogeneous singularities is free, then $\sum_{r\ge5}(r^2-5r+4)n_r-2n_2-2n_3-3$ is a perfect square.
\end{remark}

\begin{example}
For the conic-line arrangement from \cite[Example~3.2]{Po25} we get $\Delta(E_\cA)=-15$, so $\Pfrak(\cA;t)$ does not even split over the reals.
For the free arrangement from Example~\ref{ex:freeA} we get $\Delta(E_\cA)=-2\cdot6-2\cdot4+0\cdot6+3\cdot8=4=(6-4)^2$.
\end{example}

\begin{example}\label{ex:free-conics}
Consider the following six smooth conics in $\mathbb{P}^{2}_{\mathbb{C}}$:
\begin{align*}
C_0 &\colon\ xy-3xz+2yz=0, & C_3 &\colon\ 5x^2-4xy-2xz+yz=0,\\
C_1 &\colon\ 3xy-xz-2yz=0, & C_4 &\colon\ 5x^2-2xy-4xz+yz=0,\\
C_2 &\colon\ 5xy+5xz-2yz=0, & C_5 &\colon\ 5x^2-xy-xz+yz=0,
\end{align*}
and let $\mathcal{A}=\{C_0,\ldots,C_5\}$, so that $\mathcal{A}\colon f=0$ with $\deg f=12$.
The singular points of $\mathcal{A}$ are the following nine points.
\begin{center}
\begin{tabular}{lll}
\hline
point & conics through the point & type \\
\hline
$[0:1:0]$, $[0:0:1]$ & $C_0,C_1,C_2,C_3,C_4,C_5$ & ordinary $6$-fold point \\
$[1:1:1]$ & $C_0,C_1,C_3,C_4$ & ordinary $4$-fold point \\
$[3:5:-15]$ & $C_0,C_2,C_3,C_5$ & ordinary $4$-fold point \\
$[3:-15:5]$ & $C_1,C_2,C_4,C_5$ & ordinary $4$-fold point \\
$[1:0:0]$ & $C_0,C_1,C_2$ & ordinary triple point \\
$[1:3:-1]$ & $C_0,C_4,C_5$ & ordinary triple point \\
$[1:-1:3]$ & $C_1,C_3,C_5$ & ordinary triple point \\
$[1:5:5]$ & $C_2,C_3,C_4$ & ordinary triple point \\
\hline
\end{tabular}
\end{center}
Indeed, any two conics of $\mathcal{A}$ pass through exactly four of these points, so by B\'ezout's
theorem they meet transversally there and nowhere else. Hence all singular points of $\mathcal{A}$
are ordinary, $(n_3,n_4,n_6)=(4,3,2)$, $n_r=0$ otherwise, and
\[
\mu(\mathcal{A})=2\cdot 25+3\cdot 9+4\cdot 4=93 .
\]
Using \texttt{SINGULAR}~\cite{Sing} we obtain $\tau(\mathcal{A})=93$ and $\operatorname{mdr}(f)=4$.
Since $\tau_p\leq\mu_p$ for every singular point $p$, the equality $\tau(\mathcal{A})=\mu(\mathcal{A})$
shows that all singular points of $\mathcal{A}$ are quasi-homogeneous. Moreover,
\[
\gadd{(\deg f-1)^2-d_1(\deg f-1-d_1)}=121-4\cdot 7=93=\tau(\mathcal{A}),
\]
so $\mathcal{A}$ is free with exponents $(d_1,d_2)=(4,7)$ by Theorem~\gadd{\ref{thm:dPW}}. Equivalently,
$\operatorname{AR}(f)\cong S(-4)\oplus S(-7)$. Thus $\mathcal{A}$ is a free conic arrangement with
only ordinary quasi-homogeneous singularities, which, to our knowledge, is the first known example of this kind.

Let us check that this agrees with our results. Theorem~\gadd{\ref{thm:B}} with $d=2$ and $k=6$ gives
\[
\sum_{r\geq 2}(r^2-5r+4)n_r-3=2\cdot 10+4\cdot(-2)-3=9=(7-4)^2,
\]
which is the perfect square from Remark~\gadd{\ref{rem:disc} for conic arrangements}. Furthermore,
$\mathfrak{P}(\mathcal{A};t)=1+11t+28t^2=(1+4t)(1+7t)$, so Corollary~\gadd{\ref{cor:A}} gives
$e(M(\mathcal{A}))=3\cdot 6=18$, which is coherent with Proposition~\gadd{\ref{prop:comb}}:
$e(M(\mathcal{A}))=3-12+(2\cdot 5+3\cdot 3+4\cdot 2)=18$.
\end{example}
\section{The addition technique revisited}\label{sec:C}

Throughout this section $C_1,C_2\subset\PP^2_\CC$ are reduced curves without common components with $\deg C_i=c_i$.
We write $C=C_1\cup C_2$, $r=|C_1\cap C_2|$, and $i_p=i_p(C_1,C_2)$ for the intersection index of $C_1$ and $C_2$ at $p$.
Note that, according to \cite{Di24}, we have
\[
\varepsilon(C_1,C_2)=\sum_{p\in C_1\cap C_2}\bigl(\varepsilon_p(C)-\varepsilon_p(C_1)-\varepsilon_p(C_2)\bigr),
\]
where $\varepsilon_p=\mu_p-\tau_p$ and $\varepsilon_p(C_i)=0$ if $p$ is a smooth point of $C_i$.
Indeed, away from $C_1\cap C_2$ the germ of $C$ is a germ of $C_1$ or of $C_2$.

\begin{lemma}\label{lem:add}
One has $T\langle C\rangle=T\langle C_1\rangle\cap T\langle C_2\rangle$, and there is an exact sequence
\begin{equation}\label{eq:add}
0\to T\langle C\rangle\to T\langle C_1\rangle\oplus T\langle C_2\rangle\to T_{\PP^2_\CC}\to Q\to0,
\end{equation}
where the middle map is $(\theta_1,\theta_2)\mapsto\theta_1-\theta_2$ and $Q$ is a sheaf of finite length supported on the points $p\in C_1\cap C_2$ with $i_p\ge2$.
We denote the length of $Q$ by $\ell(C_1,C_2)$.
\end{lemma}

\begin{proof}
First note that the condition $\theta(h)\in(h)$ does not depend on the choice of
the local equation $h$, since $\theta(uh)=\theta(u)h+u\,\theta(h)$ for a unit $u$.
Let $p\in\PP^2_\CC$ and let $g=g_1g_2$ be a local equation of $C$ at $p$, where $g_i$
is a local equation of $C_i$ at $p$, with $g_i$ a unit if $p\notin C_i$.
Since $C_1$ and $C_2$ have no common components, $g_1$ and $g_2$ are coprime
in the factorial ring $\cO_p$.
If $\theta\in T\langle C_1\rangle\cap T\langle C_2\rangle$, then
$\theta(g)=\theta(g_1)g_2+g_1\theta(g_2)\in(g)$.
Conversely, if $\theta(g)\in(g)$, then $\theta(g_1)g_2\in(g_1)$, hence
$\theta(g_1)\in(g_1)$, and similarly $\theta(g_2)\in(g_2)$.
This proves $T\langle C\rangle=T\langle C_1\rangle\cap T\langle C_2\rangle$,
see also \cite[Section~1.4]{ST}.

Next, for subsheaves $\mathcal{A},\mathcal{B}$ of a sheaf $\mathcal{F}$, the sequence
\[
0\to\mathcal{A}\cap\mathcal{B}\to\mathcal{A}\oplus\mathcal{B}\to\mathcal{F}
\to\mathcal{F}/(\mathcal{A}+\mathcal{B})\to0,
\]
with maps $\theta\mapsto(\theta,\theta)$ and $(\theta_1,\theta_2)\mapsto\theta_1-\theta_2$,
is exact. Applying it to $\mathcal{A}=T\langle C_1\rangle$, $\mathcal{B}=T\langle C_2\rangle$
and $\mathcal{F}=T_{\PP^2_\CC}$ gives \eqref{eq:add} with
$Q=T_{\PP^2_\CC}/(T\langle C_1\rangle+T\langle C_2\rangle)$.
Since $T\langle C_i\rangle\cong E_{C_i}(1)$ is locally free by Lemma~\ref{lem:syz}
and \cite{DS}, the sheaf $Q$ is coherent.

If $p\notin C_i$, then $g_i$ is a unit, so $T\langle C_i\rangle_p=T_{\PP^2_\CC,p}$ and $Q_p=0$.
Hence $Q$ is supported on the finite set $C_1\cap C_2$, and therefore has finite length.
Finally, let $p\in C_1\cap C_2$ with $i_p=1$.
Since $i_p\ge\mult_p(C_1)\mult_p(C_2)$, with equality if and only if $C_1$ and $C_2$
have no common tangent at $p$, both curves are smooth at $p$ and meet transversally there.
Hence $dg_1$ and $dg_2$ are linearly independent at $p$, so $x=g_1$ and $y=g_2$
are local analytic coordinates centred at $p$.
A vector field $\theta=a\partial_x+b\partial_y$ satisfies $\theta(x)\in(x)$ if and only if
$a\in(x)$, so $T\langle C_1\rangle_p=\langle x\partial_x,\partial_y\rangle$, and similarly
$T\langle C_2\rangle_p=\langle\partial_x,y\partial_y\rangle$.
The sum of these two modules contains $\partial_x$ and $\partial_y$, so $Q_p=0$.
\end{proof}

\begin{proposition}\label{prop:add}
One has
\[
\ell(C_1,C_2)=\tau(C)-\tau(C_1)-\tau(C_2)-c_1c_2=c_1c_2-r-\varepsilon(C_1,C_2).
\]
In particular, $\tau(C_1\cup C_2)\ge\tau(C_1)+\tau(C_2)+c_1c_2$.
\end{proposition}

\begin{proof}
Set $K=T\langle C_1\rangle+T\langle C_2\rangle\subset T_{\PP^2_\CC}$.
By Lemma~\ref{lem:add}, the sequence \eqref{eq:add} splits into the two short exact sequences
\[
0\to T\langle C\rangle\to T\langle C_1\rangle\oplus T\langle C_2\rangle\to K\to 0
\qquad\text{and}\qquad
0\to K\to T_{\PP^2_\CC}\to Q\to 0 .
\]
Applying \ref{C1} to both of them and \ref{C2} to $Q$, we obtain
\[
c(T\langle C\rangle)\,c(T_{\PP^2_\CC})
=c(T\langle C_1\rangle)\,c(T\langle C_2\rangle)\,\bigl(1-\ell(C_1,C_2)\,[\pt]\bigr).
\]
Let $d=c_1+c_2=\deg C$. By the Euler sequence, $c(T_{\PP^2_\CC})=(1+h)^3=1+3h+3h^2$.
Next, by Proposition~\ref{prop:chern} and \eqref{eq:twist}, a reduced curve $X$ of degree $n$ satisfies
\[
c(T\langle X\rangle)=1+(3-n)h+\bigl(n^2-3n+3-\tau(X)\bigr)h^2 .
\]
The coefficients of $h$ on both sides are equal to $6-d$. Comparing the coefficients of $h^2$, we get
\[
d^2-6d+15-\tau(C)
=c_1^2+c_2^2+c_1c_2-6d+15-\tau(C_1)-\tau(C_2)-\ell(C_1,C_2),
\]
and since $d^2=c_1^2+2c_1c_2+c_2^2$, this is the first equality.

For the second equality, let $p\in C_1\cap C_2$. By \cite[Theorem~6.5.1]{Wa}, one has
\[
\mu_p(C)=\mu_p(C_1)+\mu_p(C_2)+2i_p-1,
\]
where $\mu_p(C_i)=0$ if $p$ is a smooth point of $C_i$.
If $p\in C\setminus(C_1\cap C_2)$, then the germ of $C$ at $p$ is the germ of $C_1$ or of $C_2$.
Summing over all points and using B\'ezout's theorem $\sum_{p\in C_1\cap C_2}i_p=c_1c_2$, we obtain
\[
\mu(C)=\mu(C_1)+\mu(C_2)+2c_1c_2-r,
\]
which is \cite[Lemma~3.1]{Po24}. Since $\tau=\mu-\varepsilon$, the definition of $\varepsilon(C_1,C_2)$ gives
\[
\tau(C)-\tau(C_1)-\tau(C_2)=2c_1c_2-r-\varepsilon(C_1,C_2),
\]
and substituting this into the first equality gives the second one.
Finally, the inequality follows from $\ell(C_1,C_2)\ge0$.
\end{proof}

Now we are ready to prove Theorem~\ref{thm:C}.

\begin{proof}[Proof of Theorem~\ref{thm:C}]
A direct computation gives
\[
\Pfrak(C;t)-\Pfrak(C_1;t)-\Pfrak(C_2;t)=t-1+\bigl(2c_1c_2-1-\tau(C)+\tau(C_1)+\tau(C_2)\bigr)t^2 .
\]
By Proposition~\ref{prop:add}, the coefficient of $t^2$ is equal to $c_1c_2-1-\ell(C_1,C_2)=r-1+\varepsilon(C_1,C_2)$.

Assume now that $C$ is quasi-homogeneous at $p\in C_1\cap C_2$.
By \cite{Sa71}, in suitable local coordinates a local equation of $C$ is weighted homogeneous.
Its irreducible factors are of the form $x$, $y$, or $x^a-\lambda y^b$ with $\lambda\neq0$, and all of them are weighted homogeneous with respect to the same weights.
Hence the germs of $C_1$ and $C_2$ at $p$ are quasi-homogeneous as well, so $\varepsilon_p(C)=\varepsilon_p(C_1)=\varepsilon_p(C_2)=0$.
This shows that $\varepsilon(C_1,C_2)=0$, which completes the proof.
\end{proof}

\begin{remark}
\begin{itemize}
    \item In the quasi-homogeneous case, Proposition~\ref{prop:add} tells us that $\ell(C_1,C_2)=\sum_p(i_p-1)$ is the excess intersection of $C_1$ and $C_2$.
For instance, for the line $y=0$ and the conic $y=x^2$, tangent at the origin $p$, one has $T\langle C_1\rangle+T\langle C_2\rangle=\langle\partial_x,x\partial_y,y\partial_y\rangle$ near $p$, so $\length Q_p=1=i_p-1$.

\item Recall that for a line arrangement $\cA$ all singular points are ordinary and quasi-homogeneous, and one has $\Pfrak(\cA;t)=\pi_0(\cA;t):=\pi(\cA;t)/(1+t)$, where $\pi(\cA;t)$ is the Poincar\'e polynomial of the intersection lattice of the central arrangement in $\CC^3$ defined by $\cA$.

Let $H\in\cA$ be a line, $\cA'=\cA\setminus\{H\}$, and let $\cA''$ be the restriction of the central arrangement defined by $\cA$ to the plane corresponding to $H$. The restriction consists of $r=|\cA'\cap H|$ lines in $\CC^2$.
Hence Theorem~\ref{thm:C} applied to $C_1=\cA'$ and $C_2=H$ gives us
\[
\pi_0(\cA;t)=\pi_0(\cA';t)+t\,\pi_0(\cA'';t),\qquad\text{where }\pi_0(\cA'';t)=1+(r-1)t,
\]
which is the reduced version of \cite[Theorem~1.1]{Po24} and of the classical deletion-restriction formula.
\end{itemize}
\end{remark}

Let us now focus on the Betti polynomials. Recall that if $C$ has degree $d$ and $\gadd{\kappa}$ irreducible components, then its Betti polynomial is given by
\[
B_{M(C)}(t)=1+(\gadd{\kappa}-1)t+\bigl((d-1)^2-\mu(C)-d+\gadd{\kappa}\bigr)t^2,
\]
see \cite[Theorem~3.1]{Di25} and \cite{Po24}.
This formula also follows from Proposition~\ref{prop:eC}, since $b_1(M(C))=\gadd{\kappa}-1$, and $b_3(M(C))=b_4(M(C))=0$ because $M(C)$ is affine. After comparing coefficients we obtain
\begin{equation}\label{eq:betti}
\Pfrak(C;t)=B_{M(C)}(t)+(d-\gadd{\kappa})\,t(t+1)+\varepsilon(C)\,t^2 ,
\end{equation}
which generalizes \cite[(6)]{Po25}.
In fact, \eqref{eq:betti} tells us that Theorem~\ref{thm:C} is equivalent to \cite[Theorem~1.3]{Po24} in the quasi-homogeneous case, and additionally it holds without any assumption on the singularities. Moreover, at $t=-1$ both statements reduce to the identity $e(C_1\cup C_2)=e(C_1)+e(C_2)-r$.

Finally, let us explain the role of \eqref{eq:add} in the construction of free curves.

\begin{proposition}\label{prop:del}
Assume that $C_2$ is smooth. Then there is an exact sequence
\[
0\to T\langle C_1\cup C_2\rangle\to T\langle C_1\rangle\to\cO_{C_2}(c_2)(-W)\to0,
\]
where $W$ is an effective divisor on $C_2$ of degree $\ell(C_1,C_2)$.
\end{proposition}

\begin{proof}
Let $f_2$ be an equation of $C_2$, and set $\mathcal{A}=T\langle C_1\rangle$ and
$\mathcal{B}=T\langle C_2\rangle$.
Consider the map $\nu\colon T_{\PP^2_\CC}\to\cO_{C_2}(c_2)$ given by
$\theta\mapsto\theta(f_2)|_{C_2}$, i.e., the composition
$T_{\PP^2_\CC}\to T_{\PP^2_\CC}|_{C_2}\to N_{C_2/\PP^2_\CC}\cong\cO_{C_2}(c_2)$.
It is well defined via the Euler sequence, since the Euler vector field sends $f_2$ to
$c_2f_2$, which vanishes on $C_2$.
By definition, $\ker\nu=\mathcal{B}$.
Since $C_2$ is smooth, the partial derivatives of $f_2$ have no common zero on $C_2$,
so $\nu$ is surjective and $T_{\PP^2_\CC}/\mathcal{B}\cong\cO_{C_2}(c_2)$.

By Lemma~\ref{lem:add}, the restriction $\nu|_{\mathcal{A}}$ has kernel
$\mathcal{A}\cap\mathcal{B}=T\langle C\rangle$, and its cokernel is
\[
\bigl(T_{\PP^2_\CC}/\mathcal{B}\bigr)\big/\bigl((\mathcal{A}+\mathcal{B})/\mathcal{B}\bigr)
\cong T_{\PP^2_\CC}/(\mathcal{A}+\mathcal{B})=Q .
\]
The image $\nu(\mathcal{A})$ is annihilated by $f_2$, so it is an
$\cO_{C_2}$-submodule of $\cO_{C_2}(c_2)$ whose quotient has finite length.
Since $C_2$ is a smooth curve, we get $\nu(\mathcal{A})=\cO_{C_2}(c_2)(-W)$ with
$W=\sum_p\length(Q_p)\,p$, so $W$ is effective, $Q\cong\cO_W$, and
$\deg W=\length Q=\ell(C_1,C_2)$.
\end{proof}

\begin{remark}\label{rem:dimca}
Multiplication by $f_2$ gives an inclusion $T\langle C_1\rangle(-c_2)\subset T\langle C\rangle$,
whose cokernel is a line bundle $\mathcal{F}$ on $C_2$, see \cite[Theorem~1.1]{D24}.
For a line $C_2$, this is the sequence behind Terao's addition-deletion theorem,
see \cite{Di17}, while Proposition~\ref{prop:del} describes the complementary inclusion
$T\langle C\rangle\subset T\langle C_1\rangle$.
Combining both, we get an exact sequence of vector bundles on $C_2$
\[
0\to\mathcal{F}\to T\langle C_1\rangle|_{C_2}\to\cO_{C_2}(c_2)(-W)\to0 .
\]
Since $\deg T\langle C_1\rangle|_{C_2}=\gadd{(3-c_1)c_2}$ and $\deg W=c_1c_2-r-\varepsilon(C_1,C_2)$
by Proposition~\ref{prop:add}, we get
$\deg\mathcal{F}=3c_2-c_2^2-r-\varepsilon(C_1,C_2)$,
coherently with \cite[Theorem~1.1]{D24}.
\end{remark}

\begin{corollary}\label{cor:addfree}
Let $C_1$ be free with exponents $(a_1,a_2)$, let $C_2$ be smooth of degree $c$, and assume that $C_1\cup C_2$ is quasi-homogeneous at all points of $C_1\cap C_2$.
If $C_1\cup C_2$ is free with exponents $(b_1,b_2)$, then
\[
b_1+b_2=a_1+a_2+c\qquad\text{and}\qquad b_1b_2=a_1a_2+(c-1)^2+r-1 .
\]
In particular, the exponents $(a_1,a_2+c)$ can occur only if $r=c(a_1-c+2)$.
For $c=1$ this is the condition $r=a_1+1$ from Terao's addition theorem, and for $c=2$ it reads $r=2a_1$.
\end{corollary}

\begin{proof}
Since $\Pfrak(C_2;t)=1+(c-1)t+(c-1)^2t^2$, the claim follows from Theorem~\ref{thm:C} by comparing the coefficients of $(1+b_1t)(1+b_2t)$ and $(1+a_1t)(1+a_2t)+\Pfrak(C_2;t)+t-1+(r-1)t^2$.
\end{proof}

\begin{example}\label{ex:addition}
Consider the arrangement $C_1$ of $7$ lines and the conic $C_2\colon x^2+y^2-2z^2=0$
from \cite[Example~3.4]{Po24}, and let $C=C_1\cup C_2$.
The curve $C_1$ is free with exponents $(3,3)$, it has six triple points and three double
points, and $\tau(C_1)=27$.
The conic passes transversally through the four triple points $[\pm1:\pm1:1]$ of $C_1$,
where $i_p=3$, and it meets the line $x=0$ transversally at two further points, so $r=6$.
All singular points of $C$ are ordinary of multiplicity at most $4$, hence quasi-homogeneous,
so $\varepsilon(C_1,C_2)=0$, and Proposition~\ref{prop:add} gives
\[
\ell(C_1,C_2)=c_1c_2-r=14-6=\sum_{p\in C_1\cap C_2}(i_p-1)
\qquad\text{and}\qquad
\tau(C)= 49,
\]
which is coherent with the weak-combinatorics $(n_2,n_3,n_4)=(5,2,4)$.
Moreover, Theorem~\ref{thm:C} gives $\Pfrak(C;t)=1+8t+15t^2=(1+3t)(1+5t)$, and
$r=6=2a_1$, which agrees with Corollary~\ref{cor:addfree} and the exponents $(3,5)$ of $C$.

Next, consider \cite[Example~3.5]{Po24}, where $C_1\colon x^3+y^3=0$ is a pencil of three
lines through $[0:0:1]$ and $C_2\colon x^3+y^3+z^3=0$ is the Fermat cubic.
Let $\omega$ be a primitive cube root of unity.
The lines of $C_1$ meet $C_2$ only at the points $p_j=[-\omega^j:1:0]$, $j=0,1,2$, which are
flexes of $C_2$, and the lines are the inflectional tangents at these points.
Hence $r=3$ and $i_{p_j}=3$ for all $j$.
The singular points of $C=C_1\cup C_2$ are the ordinary triple point $[0:0:1]$ and three
singular points of type $A_5$ at $p_0,p_1,p_2$, all of them quasi-homogeneous, so
$\varepsilon(C_1,C_2)=0$ and Proposition~\ref{prop:add} gives
\[
\ell(C_1,C_2)=9-3=6
\qquad\text{and}\qquad
\tau(C)=19=4+3\cdot5 .
\]
In order to determine the divisor $W$ in Proposition~\ref{prop:del}, note that near each
$p_j$ there are local analytic coordinates in which the line is $y=0$ and the cubic is $y=x^3$.
Then $T\langle C_1\rangle_{p_j}=\langle\partial_x,\,y\partial_y\rangle$ and
$T\langle C_2\rangle_{p_j}=\langle\partial_x+3x^2\partial_y,\,(y-x^3)\partial_y\rangle$, so
\[
T\langle C_1\rangle_{p_j}+T\langle C_2\rangle_{p_j}=\langle\partial_x,\,x^2\partial_y,\,y\partial_y\rangle,
\]
and $Q_{p_j}\cong\cO_{p_j}/(x^2,y)$ has length $2$.
Hence $W=2(p_0+p_1+p_2)$, and $\deg W=6=\ell(C_1,C_2)$, as it should be.
\end{example}

\section{Bernstein--Sato polynomials}\label{sec:BS}

In this section we look at the above results from the viewpoint of Bernstein--Sato polynomials.
The singularities that play the main role in Sections~\ref{sec:A}--\ref{sec:C}, namely quasi-homogeneous and ordinary singularities, are also singularities for which Bernstein--Sato polynomials can be computed, or at least efficiently controlled, by hand.
In the quasi-homogeneous case this is the theorem of Brian\c{c}on, Granger and Maisonobe \cite{BGM}.
In the case of ordinary singularities a single blow-up suffices, namely the one that we \gadd{use in the proof of Theorem~\ref{thm:D}}.
Our main reference on this topic is the beautiful notes by Popa \cite{PoBS}.

\subsection*{Basic notions}
Let $U$ be a complex manifold and let $h\in\cO_U(U)$ be non-invertible.
The Bernstein--Sato polynomial $b_h(s)$ is the \gdel{monique} \gadd{monic} polynomial of minimal degree such that
\[
b_h(s)\,h^s=P(s)\,h^{s+1}\qquad\text{for some }P(s)\in\sD_U[s],
\]
where $\sD_U$ denotes the sheaf of differential operators on $U$.
The local Bernstein--Sato polynomial $b_{h,x}(s)$ is equal to $b_{h|V}(s)$ for every sufficiently small neighbourhood $V$ of $x\in U$, see \cite[Definitions~1.3 and~1.5]{PoBS}.
It does not change under analytic changes of coordinates and under multiplication of $h$ by a unit \cite[Section~1]{PoBS}.

\begin{definition}\label{def:BS}
Let $C\subset\PP^2_\CC$ be a reduced curve and $p\in C$.
The local Bernstein--Sato polynomial of $C$ at $p$ is $b_{C,p}(s):=b_{h,p}(s)$, where $h$ is a local equation of $C$ at $p$, and the Bernstein--Sato polynomial of $C$ is
\[
b_C(s):=\lcm_{p\in C}\,b_{C,p}(s),
\]
see \cite[Remark~1.8]{PoBS}.
We put $\tb_{C,p}(s)=b_{C,p}(s)/(s+1)$, and
\[
R_p(C)=\{\alpha\colon b_{C,p}(-\alpha)=0\},\qquad
\tR_p(C)=\{\alpha\colon \tb_{C,p}(-\alpha)=0\},\qquad
R(C)=\bigcup_{p\in C}R_p(C).
\]
\end{definition}

By Kashiwara's theorem \cite{Ka}, see also \cite[Theorem~3.1]{PoBS}, we have $R_p(C)\subset\QQ_{>0}$.
Moreover, $s+1$ divides $b_{C,p}(s)$, and $b_{C,p}(s)=s+1$ if $p$ is a smooth point of $C$ \cite[Remark~1.9]{PoBS}.
Hence $b_C(s)$ is the least common multiple of $s+1$ and of the finitely many polynomials $b_{C,p}(s)$ with $p\in\Sing(C)$, and $R(C)$ is the set of roots of $b_C(-s)$.
The log canonical threshold and the minimal exponent of $C$ at $p$ are then given as
\[
\lct_p(C)=\min R_p(C)\qquad\text{and}\qquad\talpha_p(C)=\min\tR_p(C),
\]
with the convention $\min\emptyset=+\infty$.
The first equality is the theorem of Yano, Lichtin and Koll\'ar \cite{Ko}, see also \cite[Theorem~7.2]{PoBS}, and the second one is the definition that was taken from \cite[Section~9]{PoBS}.
Recall that 
$$\lct_p(C)=\min\{1,\talpha_p(C)\},$$ and that $$\lct(C)=\min_{p\in C}\lct_p(C)=\min R(C).$$

We will use the following four facts, valid for every singular point $p$ of $C$.
\begin{enumerate}[label=\textup{(B\arabic*)},leftmargin=*]
\item\label{B1} Every jumping number of $C$ at $p$ lying in $(0,1]$ belongs to $R_p(C)$, see \cite{ELSV} and \cite[Theorem~7.3]{PoBS}.
\item\label{B2} Let $U$ be a small neighbourhood of $p$ and let $\sigma\colon V\to U$ be a log resolution of $C\cap U$ which is an isomorphism over $U\setminus C$.
Write $\sigma^*C=\sum_ia_iE_i$ and $K_{V/U}=\sum_ib_iE_i$.
Then every element of $R_p(C)$ is of the form $(b_i+1+\ell)/a_i$ for some $i$ and some integer $\ell\ge0$, see \cite[Theorem~5]{Lich} and \cite[Theorem~4.4]{PoBS}.
\item\label{B3} One has $\tR_p(C)\subset\bigl[\talpha_p(C),\,2-\talpha_p(C)\bigr]$, see \cite[Theorem~0.4]{SaM} and \cite[Remark~9.11]{PoBS}.
\item\label{B4} If $m=\mult_p(C)$, then $\talpha_p(C)\le 2/m$, with equality if $p$ is an ordinary singular point, see \cite[Theorem~E]{MP} and \cite[Theorem~9.9(3) and Example~9.10]{PoBS}.
\end{enumerate}
Since $m\ge2$ at a singular point, \ref{B4} gives $\talpha_p(C)\le1$, and therefore
\begin{equation}\label{eq:lctalpha}
\lct_p(C)=\talpha_p(C)\le\frac{2}{\mult_p(C)}\qquad\text{for every }p\in\Sing(C).
\end{equation}
In particular, for singular points of reduced plane curves the minimal exponent is equal to the log canonical threshold, and the additional information is contained in the remaining roots.

\subsection*{Quasi-homogeneous singularities}
Let $p\in\Sing(C)$ be a quasi-homogeneous singular point.
By \cite{Sa71} there are local analytic coordinates $(x,y)$ centred at $p$ and weights $w=(w_1,w_2)\in\QQ^2_{>0}$ such that $C$ is defined near $p$ by a polynomial $h$ which is weighted homogeneous of degree $1$, i.e., $w_1\alpha+w_2\beta=1$ for every monomial $x^\alpha y^\beta$ occurring in $h$.
We put $|w|=w_1+w_2$, and we denote by $\Sigma_p$ the set of weights of the monomials in a monomial basis of the Milnor algebra $M_p=\CC\{x,y\}/\langle\partial_xh,\partial_yh\rangle$.
The result below gathers some elementary properties that are well-known to the experts, but we present them for the completness of the exposition.
\begin{proposition}\label{prop:qhBS}
In the above setting, the following hold.
\begin{enumerate}[label=\textup{(\alph*)},leftmargin=*]
\item $b_{C,p}(s)=(s+1)\prod_{\rho\in\Sigma_p}(s+|w|+\rho)$. In particular, all roots of $\tb_{C,p}(s)$ are simple.
\item $\lct_p(C)=\talpha_p(C)=|w|$.
\item The set $\tR_p(C)$ is invariant under $\alpha\mapsto2-\alpha$. In particular, $\max\tR_p(C)=2-\lct_p(C)$.
\item $\mu_p=\tau_p=(1-w_1)(1-w_2)/(w_1w_2)$.
\end{enumerate}
\end{proposition}

\begin{proof}
Part (a) is \cite{BGM}, see \cite[Theorem~2.4]{PoBS}.
Since $\rho(1)=0$ is the smallest element of $\Sigma_p$, the greatest root of $\tb_{C,p}(s)$ is $-|w|$, see \cite[Example~9.4]{PoBS}, so (b) follows from \eqref{eq:lctalpha}.

For (c) and (d), note that $\partial_xh$ and $\partial_yh$ are weighted homogeneous of degrees $1-w_1$ and $1-w_2$.
Their only common zero is the origin.
Indeed, by the Euler relation $h=w_1x\,\partial_xh+w_2y\,\partial_yh$ their common zeros are singular points of $C$, and their common zero locus is invariant under the $\CC^*$-action defined by $w$, whose orbits accumulate at the isolated singular point $p$.
Hence $\partial_xh,\partial_yh$ form a regular sequence in $\CC[x,y]$, the graded algebra $A=\CC[x,y]/\langle\partial_xh,\partial_yh\rangle$ is finite dimensional and isomorphic to $M_p$, and its Hilbert series is
\[
H_p(T)=\sum_\rho\dim A_\rho\,T^\rho=\prod_{i=1}^2\frac{1-T^{1-w_i}}{1-T^{w_i}} .
\]
In particular, $\Sigma_p=\{\rho\colon A_\rho\neq0\}$.
A direct computation gives $T^{2-2|w|}H_p(T^{-1})=H_p(T)$, so $\Sigma_p$ is invariant under $\rho\mapsto2-2|w|-\rho$.
By (a) we have $\tR_p(C)=\{|w|+\rho\colon\rho\in\Sigma_p\}$, which gives (c).
Finally, $\mu_p=\dim A=H_p(1)=\prod_i(1-w_i)/w_i$\gadd{, see also \cite{MO}}, and $\mu_p=\tau_p$ by \cite{Sa71}.
\end{proof}

As a consequence, we obtain the following a bit technical description of curves with $e(M(C))=\Pfrak(C;-1)$.

\begin{corollary}\label{cor:ABS}
Let $C\subset\PP^2_\CC$ be a reduced curve of degree $d$ such that $e(M(C))=\Pfrak(C;-1)$, i.e., $C$ admits only quasi-homogeneous singularities.
For $p\in\Sing(C)$ let $w_p=(w_{1,p},w_{2,p})$ be the weights at $p$. Then the following hold.
\begin{enumerate}[label=\textup{(\alph*)},leftmargin=*]
\item $b_C(s)$ is the least common multiple of $s+1$ and of the polynomials $(s+1)\prod_{\rho\in\Sigma_p}(s+|w_p|+\rho)$ with $p\in\Sing(C)$, and every root of $b_C(s)$ other than $-1$ is simple.
\item For every $p\in\Sing(C)$ the set $\tR_p(C)$ is symmetric with respect to $1$. In particular, $\tb_{C,p}\bigl(\lct_p(C)-2\bigr)=0$.
\item One has
\[
e(M(C))=d^2-3d+3-\sum_{p\in\Sing(C)}\frac{(1-w_{1,p})(1-w_{2,p})}{w_{1,p}\,w_{2,p}} .
\]
\item If $C$ is free with exponents $(d_1,d_2)$, then
\[
\sum_{p\in\Sing(C)}\frac{(1-w_{1,p})(1-w_{2,p})}{w_{1,p}\,w_{2,p}}=(d-1)^2-d_1d_2 .
\]
\end{enumerate}
\end{corollary}

\begin{proof}
Parts (a) and (b) follow from Proposition~\ref{prop:qhBS}, since the multiplicity of a root of a least common multiple is the maximum of its multiplicities, and $-1$ is a root of $b_{C,p}(s)$ of multiplicity at most two.
Part (c) follows from Theorem~\ref{thm:A}, \eqref{eq:twist} and Proposition~\ref{prop:qhBS}(d), and (d) follows from (c) and Corollary~\ref{cor:A}.
\end{proof}

\begin{example}
Let $C\colon y^2z-x^3=0$ be the cuspidal cubic.
Its only singular point is a cusp with weights $w=(1/3,1/2)$ and $\Sigma_p=\{0,1/3\}$, so
\[
b_C(s)=(s+1)\bigl(s+\tfrac56\bigr)\bigl(s+\tfrac76\bigr).
\]
Observe that $\tR_p(C)=\{5/6,7/6\}$ is symmetric with respect to $1$, and Corollary~\ref{cor:ABS}(c) gives $e(M(C))=3-2=1$.
\end{example}

Theorem~\ref{thm:A} characterizes the curves with $e(M(C))=\Pfrak(C;-1)$ by the vanishing of $\mu(C)-\tau(C)$.
Corollary~\ref{cor:ABS}(b) leads us to the following question.

\begin{question}\label{q:sym}
Let $p$ be a singular point of a reduced plane curve $C$ such that $2-\talpha_p(C)\in\tR_p(C)$.
Is $p$ quasi-homogeneous?
\end{question}
If the answer is positive, then Theorem~\ref{thm:A} and Corollary~\ref{cor:ABS} give
\[
e(M(C))=\Pfrak(C;-1)\iff\tb_{C,p}\bigl(\lct_p(C)-2\bigr)=0\ \text{ for every }p\in\Sing(C),
\]
i.e., the correction term $\mu(C)-\tau(C)$ in Theorem~\ref{thm:A} vanishes if and only if every set $\tR_p(C)$ attains the upper end point of the interval in \ref{B3}.

\subsection*{Ordinary singularities and the proof of Theorem~\ref{thm:D}}
For ordinary singularities the roots can be located without any assumption on quasi-homogeneity.
We need the following well-known characterization of quasi-homogeneous ordinary singular points, and for the completeness of our presentation we deliver a short proof.
\begin{lemma}\label{lem:ordqh}
Let $p$ be an ordinary $r$-fold point of $C$.
Then $p$ is quasi-homogeneous if and only if, after a suitable change of local analytic coordinates, $C$ is defined near $p$ by a homogeneous polynomial of degree $r$.
\end{lemma}

\begin{proof}
A homogeneous polynomial of degree $r$ is weighted homogeneous with $w=(1/r,1/r)$.
Conversely, assume that $p$ is quasi-homogeneous, and let $h$ be a local equation of $C$ which is weighted homogeneous of degree $1$ with respect to $w$, as above.
For $r=2$ the point $p$ is a node, which is given by $xy=0$ in suitable coordinates.
Let $r\ge3$ and let $h_r$ be the homogeneous part of $h$ of degree $r$.
Since $p$ is an ordinary singular point, $h_r$ is a product of $r$ pairwise distinct linear forms; this property does not depend on the choice of coordinates, since $h_r$ transforms via the linear part of a change of coordinates.
Every monomial $x^\alpha y^\beta$ of $h_r$ satisfies $\alpha+\beta=r$ and $w_1\alpha+w_2\beta=1$.
If $w_1\neq w_2$, then $h_r$ is a monomial, which is impossible for $r\ge3$.
Hence $w_1=w_2=1/r$ and $h=h_r$ is homogeneous of degree $r$.
\end{proof}

\begin{proof}[Proof of Theorem~\ref{thm:D}]
By \ref{B4} and \eqref{eq:lctalpha}, we have $\lct_p(C)=\talpha_p(C)=2/r$.

Let $U$ be a small neighbourhood of $p$ and let $\sigma\colon V\to U$ be the blow-up of $p$
with exceptional curve $E$.
Since the branches of $C$ at $p$ are smooth and pairwise transverse, $\sigma$ is a log
resolution of $C\cap U$ which is an isomorphism over $U\setminus\{p\}$, with
\[
\sigma^*C=\widetilde C+rE\qquad\text{and}\qquad K_{V/U}=E,
\]
where the strict transform $\widetilde C$ is smooth and meets $E$ transversally.
Denote by $\cJ(c\,C)$ the multiplier ideal of $c\,C$.
For $0<c<1$ we get
\[
\cJ(c\,C)=\sigma_*\cO_V\bigl(K_{V/U}-\lfloor c\,\sigma^*C\rfloor\bigr)
=\sigma_*\cO_V\bigl((1-\lfloor cr\rfloor)E\bigr)
=\mathfrak m_p^{\max\{0,\lfloor cr\rfloor-1\}},
\]
while $\cJ(C)=\cO_U(-C)\subsetneq\mathfrak m_p^{r-2}$, see \cite[Section~6]{PoBS}.
Hence the jumping numbers of $C$ at $p$ in $(0,1]$ are exactly the numbers $j/r$ with
$2\le j\le r$, and the first inclusion \gadd{for $j\le r$} follows from \ref{B1}.

For the second inclusion, \ref{B2} applied to $\sigma$ shows that every element of $R_p(C)$
is either of the form $(2+\ell)/r$ with $\ell\ge0$, coming from $E$, for which $a_i=r$ and
$b_i=1$, or a positive integer, coming from the components of $\widetilde C$, for which
$a_i=1$ and $b_i=0$.
By \ref{B3}, we have $\tR_p(C)\subset[2/r,\,2-2/r]$, and the only integer in this interval
is $1=r/r$.
Since $R_p(C)=\tR_p(C)\cup\{1\}$, the second inclusion follows.

It remains to show that $(r+1)/r\in R_p(C)$ for $r\ge3$.
By results of Malgrange \gadd{\cite{Ma}} and Varchenko, see \cite[(1)]{SaBP}, there are non-negative integers $\nu_k$ such that $\tR_p(C)=\{\alpha_k-\nu_k\}_k$ and $\min_k(\alpha_k-\nu_k)=\min_k\alpha_k$, where the $\alpha_k$ are the spectral numbers of $C$ at $p$.
If $h_r$ is the homogeneous part of degree $r$ of a local equation $h$, then $h_r+t(h-h_r)$, $t\in[0,1]$, is a $\mu$-constant family of ordinary $r$-fold points, so the spectral numbers of $C$ at $p$ are those of $h_r$, namely $j/r$ with $2\le j\le2r-2$, see \cite[Introduction and (1.1.20)]{SaBP}.
Since $(r+1)/r-1<2/r=\min_k\alpha_k$, the spectral number $(r+1)/r$ is not shifted, i.e., $(r+1)/r\in\tR_p(C)$.

Assume now that $p$ is quasi-homogeneous.
By Lemma~\ref{lem:ordqh}, the curve $C$ is defined near $p$ by a homogeneous polynomial $h$
of degree $r$ with an isolated singularity, i.e., the germ of $C$ at $p$ is analytically
equivalent to a central arrangement of $r$ pairwise distinct lines in $\CC^2$.
Since any two of these lines are linearly independent, this arrangement is generic, and the
formula for $b_{C,p}(s)$ is the case $n=2$ of Walther's formula for generic central
arrangements in $\CC^n$ \cite{Wal}, see also \cite{SaH}.
For the convenience of the reader, we deduce it from Proposition~\ref{prop:qhBS}(a).
We have $w=(1/r,1/r)$ and $|w|=2/r$.
The partial derivatives $\partial_xh,\partial_yh$ form a regular sequence of forms of
degree $r-1$, so the Milnor algebra has Hilbert series $(1+T+\dots+T^{r-2})^2$ with respect
to the standard grading, which has a nonzero coefficient in every degree $0\le k\le2r-4$.
Thus $\Sigma_p=\{k/r\colon0\le k\le2r-4\}$, and Proposition~\ref{prop:qhBS}(a) gives
\[
b_{C,p}(s)=(s+1)\prod_{k=0}^{2r-4}\Bigl(s+\frac{k+2}r\Bigr)
=(s+1)\prod_{j=2}^{2r-2}\Bigl(s+\frac jr\Bigr),
\]
which completes the proof.
\end{proof}

\begin{remark}\label{rem:ordBS}
Since ordinary singular points of multiplicity $r\le4$ are quasi-homogeneous, and so are all singular points of line arrangements, Theorem~\ref{thm:D} gives
\[
(s+1)^2,\qquad (s+1)^2\bigl(s+\tfrac23\bigr)\bigl(s+\tfrac43\bigr),\qquad (s+1)^2\bigl(s+\tfrac12\bigr)\bigl(s+\tfrac34\bigr)\bigl(s+\tfrac54\bigr)\bigl(s+\tfrac32\bigr)
\]
for $r=2,3,4$, respectively.
For $r\ge5$ both inclusions in Theorem~\ref{thm:D} can be strict.
For the ordinary $5$-fold point of the sextic from Example~\ref{ex:nonqh}, \texttt{SINGULAR} computations give $R_p(C)=\{j/5\colon2\le j\le7\}$, and for the ordinary $6$-fold point $x^6+y^6+x^4y^4=0$, with $\mu_p=25$ and $\tau_p=24$, it gives $R_p=\{j/6\colon2\le j\le9\}$.
In both cases the only missing candidate is $(2r-2)/r=2-\talpha_p$, compare it with Question~\ref{q:sym}.

Moreover, $R_p(C)$ need not consist of consecutive multiples of $1/r$: for the ordinary $10$-fold point $\frac1{10}x^{10}+\frac1{10}y^{10}+x^3y^8+x^8y^3+8x^6y^6-128x^7y^7=0$ one has $\gadd{R_p}=\{j/10\colon2\le j\le14\}\cup\{16/10\}$, see \cite[page 7]{SaBP}.
\end{remark}

Taking the least common multiple of the local Bernstein--Sato polynomials, we obtain the following.

\begin{corollary}\label{cor:arrBS}
Let $C\subset\PP^2_\CC$ be a reduced curve whose singular points are all ordinary, assume that $\Sing(C)\neq\emptyset$, and let $m$ be the maximal multiplicity of a singular point of $C$. Then $\lct(C)=2/m$ and
\[
\bigcup_{n_r>0}\Bigl\{\frac jr\colon 2\le j\le \min\{r+1,2r-2\}\Bigr\}\subseteq R(C)\subseteq\bigcup_{n_r>0}\Bigl\{\frac jr\colon 2\le j\le 2r-2\Bigr\}.
\]
Moreover, if all singular points of $C$ are quasi-homogeneous, for instance if $m\le4$ or if $C$ is a line arrangement, then the second inclusion is an equality and
\[
b_C(s)=(s+1)^2\prod_{\alpha\in R(C)\setminus\{1\}}(s+\alpha).
\]
\end{corollary}

\begin{proof}
Since $b_C(s)$ is the least common multiple of the local Bernstein--Sato polynomials, this follows from Theorem~\ref{thm:D} and Remark~\ref{rem:ordBS}.
\end{proof}

Observe that in the quasi-homogeneous case $b_C(s)$ depends only on the set $\{r\colon n_r>0\}$, whereas $\Pfrak(C;t)=\Pcomb(C;t)$ depends on the numbers $n_r$ themselves, see Corollary~\ref{cor:comb}.

\begin{example}\label{ex:BSarr}
For the free arrangement $\cA$ of $9$ lines and one conic from Example~\ref{ex:freeA} we have $(n_2,n_3,n_4)=(6,4,6)$, so Corollary~\ref{cor:arrBS} gives
\[
b_\cA(s)=(s+1)^2\bigl(s+\tfrac12\bigr)\bigl(s+\tfrac23\bigr)\bigl(s+\tfrac34\bigr)\bigl(s+\tfrac54\bigr)\bigl(s+\tfrac43\bigr)\bigl(s+\tfrac32\bigr)
\]
and $\lct(\cA)=1/2$.
In Example~\ref{ex:addition}, the arrangement $C_1$ of $7$ lines has six triple points and three double points, while $C=C_1\cup C_2$ has $(n_2,n_3,n_4)=(5,2,4)$.
Hence
\[
b_{C_1}(s)=(s+1)^2\bigl(s+\tfrac23\bigr)\bigl(s+\tfrac43\bigr),
\]
and $b_C(s)$ is equal to $b_\cA(s)$ above.
\end{example}

\subsection*{Free \texorpdfstring{$d$}{d}-arrangements and the proof of Corollary~\ref{cor:D}}

\begin{proof}[Proof of Corollary~\ref{cor:D}]
Since $r^2-5r+4=(r-1)(r-4)$ is equal to $-2$ for $r\in\{2,3\}$ and to $0$ for $r=4$, Theorem~\ref{thm:B} can be written as
\begin{equation}\label{eq:BBS}
\sum_{r\ge5}(r-1)(r-4)\,n_r=3+3dk(d-2)+2(n_2+n_3)+(d_2-d_1)^2 .
\end{equation}
The right-hand side of \eqref{eq:BBS} is at least $3$, so $n_r>0$ for some $r\ge5$, and hence $m\ge5$.
Each curve $C_i$ is smooth, so it contributes exactly one branch at every point of $\cA$ it passes through.
Thus an $r$-fold point of $\cA$ lies on exactly $r$ curves of $\cA$, and $m\le k$.
By Corollary~\ref{cor:arrBS}, we have $\lct(\cA)=2/m\le2/5$.
Finally, since all singular points of $\cA$ are quasi-homogeneous, Theorem~\ref{thm:D} gives $R(\cA)\supseteq\{j/m\colon2\le j\le2m-2\}$, so $b_\cA(s)$ has at least $2m-3\ge7$ pairwise distinct roots, which completes the proof.
\end{proof}
\begin{example}
Let us come back to Example \ref{ex:free-conics}. We have
\[
b_{\mathcal{A}}(s)=(s+1)^2\Bigl(s+\tfrac13\Bigr)\Bigl(s+\tfrac12\Bigr)\Bigl(s+\tfrac23\Bigr)
\Bigl(s+\tfrac34\Bigr)\Bigl(s+\tfrac56\Bigr)\Bigl(s+\tfrac76\Bigr)\Bigl(s+\tfrac54\Bigr)
\Bigl(s+\tfrac43\Bigr)\Bigl(s+\tfrac32\Bigr)\Bigl(s+\tfrac53\Bigr),
\]
which has $11\geq 2m-3=9$ pairwise distinct roots.
\end{example}
Note that a $d$-arrangement whose only singularities are nodes satisfies $\lct(\cA)=1$.
Corollary~\ref{cor:D} tells us that the freeness of a $d$-arrangement with quasi-homogeneous singularities forces the log canonical threshold to drop to \textbf{at most} $2/5$, see also \cite[Theorem 5.2]{Po24a} and \cite[Theorem A]{Po23}.

\subsection*{The addition technique and log canonical thresholds}
The proof of Theorem~\ref{thm:C} is based on the observation that the branches of a quasi-homogeneous germ are quasi-homogeneous with respect to the same weights.
In the language of Bernstein--Sato polynomials, this leads to the following additivity property. This result looks very natural and we believe that it has to be known among the experts. However, we were unable to locate this result in the literature, and therefore we decided to present it with a complete proof.

\begin{proposition}\label{prop:addBS}
Let $C_1,C_2\subset\PP^2_\CC$ be reduced curves without common components, and let $p\in\Sing(C_1)\cap\Sing(C_2)$.
If $C=C_1\cup C_2$ is quasi-homogeneous at $p$, then $C_1$ and $C_2$ are quasi-homogeneous at $p$, and
\[
\frac1{\lct_p(C_1\cup C_2)}=\frac1{\lct_p(C_1)}+\frac1{\lct_p(C_2)} .
\]
\end{proposition}

\begin{proof}
As in the proof of Theorem~\ref{thm:C}, we choose local coordinates in which a local equation $h$ of $C$ at $p$ is weighted homogeneous of degree $1$ with respect to $w$.
Its irreducible factors are weighted homogeneous with respect to $w$, so $h=h_1h_2$, where $h_i$ is a local equation of $C_i$ which is weighted homogeneous of some degree $\rho_i>0$ with respect to $w$, and $\rho_1+\rho_2=1$.
Then $h_i$ is weighted homogeneous of degree $1$ with respect to $w/\rho_i$, and $p$ is an isolated singular point of $C_i$.
By Proposition~\ref{prop:qhBS}(b), we get $\lct_p(C)=|w|$ and $\lct_p(C_i)=|w|/\rho_i$, which completes the proof.
\end{proof}

\begin{remark}\label{rem:addBS}
\begin{itemize}
\item Consider the union of the two cusps $(y^2-x^3)(y^2-2x^3)=0$. The log canonical threshold is $5/12$, and indeed we have $(6/5+6/5)^{-1}=5/12$.
\item If $p$ is an ordinary singular point of $C$, then \eqref{eq:lctalpha} and \ref{B4} show that the formula holds whenever $C_1$ and $C_2$ are both singular at $p$, even if $C$ is not quasi-homogeneous at $p$, since $2/(r_1+r_2)=(r_1/2+r_2/2)^{-1}$.
In general, however, the quasi-homogeneity assumption cannot be dropped.
For the union of the two transverse cusps $(x^2-y^3)(y^2-x^3)=0$, which is not quasi-homogeneous at the origin, \texttt{SINGULAR} gives $\lct_0=1/2$, whereas $(6/5+6/5)^{-1}=5/12$.
This is analogous to the role played by the correction term $\varepsilon(C_1,C_2)$ in Theorem~\ref{thm:C}.
\end{itemize}
\end{remark}

\section*{Acknowledgments}
The idea of looking at the Zeuthen--Segre formula was suggested to me by Fabrizio Catanese during a mini-workshop in Hannover in March 2025 and I want to thank him for this pieace of advise. I want to thank Alex Dimca for useful comments, and Mircea Musta\c{t}\u{a} for pointing out that $-(r+1)/r$ is \gdel{always} a root in Theorem~\ref{thm:D} \gadd{for every $r\ge3$} and explaining to me the proof of that.

\section*{AI Declaration}
Example \ref{ex:free-conics} was found by symbolic experiments in positive characteristic using \texttt{Claude}, and the lifting and all checks in characteristic zero were performed by the author using \texttt{SINGULAR}.

\section*{Funding}
Piotr Pokora is supported by the National Science Centre (Poland) Sonata Bis Grant \textbf{\path{2023/50/E/ST1/00025}}. For the purpose of Open Access, the author has applied a CC-BY public copyright license to any Author Accepted Manuscript (AAM) version arising from this submission.

\noindent
Piotr Pokora, 
Department of Mathematics,
University of the National Education Commission Krakow,
Podchora\.zych 2,
PL-30-084 Krak\'ow, Poland. \\
\nopagebreak
\noindent
Email: \texttt{piotr.pokora@uken.krakow.pl}

\appendix
\section*{Appendix}
\subsection*{\texttt{SINGULAR} script for Example \ref{ex:nonqh}\texorpdfstring{{\gblock{} and Remark~\ref{rem:ordBS}}}{ and Remark}}
We verify the local invariants of the curve
$C\colon f = x^5z + y^5z + x^3y^3 = 0$ at the point $p=[0:0:1]$ using
\texttt{SINGULAR} \cite{Sing}. In the affine chart $z=1$, the point $p$ is the origin,
and we work in the local ring $\CC\{x,y\}$ via the local degree ordering
\texttt{ds}. The Milnor and Tjurina numbers are computed as
\[
  \mu(C,p)=\dim_\CC \CC\{x,y\}/(g_x,g_y),
  \qquad
  \tau(C,p)=\dim_\CC \CC\{x,y\}/(g,g_x,g_y),
\]
where $g=f(x,y,1)$.

\begin{verbatim}
ring r = 0,(x,y),ds;              // local ordering at the origin
poly f = x5 + y5 + x3y3;          // affine chart z = 1 of C
int mu  = vdim(std(jacob(f)));                // Milnor number
int tau = vdim(std(jacob(f) + ideal(f)));     // Tjurina number
printf("mu = %s, tau = %s, mu - tau = %s", mu, tau, mu - tau);
if (mu == tau) { "quasi-homogeneous (Saito)"; }
else           { "not quasi-homogeneous (Saito)"; }
\end{verbatim}

The output is $\mu(C,p)=16$ and $\tau(C,p)=15$. Hence the singularity
of $C$ at $p$ is not quasi-homogeneous, by Saito's theorem.

\gadd{Appending the following lines to the script above computes the local Bernstein--Sato polynomial used in Remark~\ref{rem:ordBS}.}
{\gblock
\begin{verbatim}
LIB "gmssing.lib";                // local Bernstein--Sato polynomial
list L = bernstein(f);            // roots of b_{f,0}(s) and multiplicities
L;
\end{verbatim}}
\gadd{The output is $-\tfrac75,-\tfrac65,-1,-\tfrac45,-\tfrac35,-\tfrac25$, where $-1$ has multiplicity two and all other roots are simple, i.e., $R_p(C)=\{j/5\colon2\le j\le7\}$.}

\subsection*{\texttt{SINGULAR} script for Remark \ref{rem:ordBS}}
We compute the local Bernstein--Sato polynomial of the ordinary $6$-fold point
$x^6+y^6+x^4y^4=0$ at the origin using the procedure \texttt{bernstein} from the library
\texttt{gmssing.lib} of \texttt{SINGULAR} \cite{Sing}, which computes the roots of the local
Bernstein--Sato polynomial of an isolated hypersurface singularity, together with their
multiplicities, via the Gauss--Manin connection.

\begin{verbatim}
LIB "gmssing.lib";                // local Bernstein--Sato polynomial
ring r = 0,(x,y),ds;              // local ordering at the origin
poly f = x6 + y6 + x4y4;          // ordinary 6-fold point
int mu  = vdim(std(jacob(f)));                // Milnor number
int tau = vdim(std(jacob(f) + ideal(f)));     // Tjurina number
printf("mu = %s, tau = %s", mu, tau);
list L = bernstein(f);            // roots of b_{f,0}(s) and multiplicities
L;
\end{verbatim}

The output is $\mu_p=25$, $\tau_p=24$, and the roots
$-\tfrac32,-\tfrac43,-\tfrac76,-1,-\tfrac56,-\tfrac23,-\tfrac12,-\tfrac13$,
where $-1$ has multiplicity two and all other roots are simple.

\subsection*{\texttt{SINGULAR} script for Remark \ref{rem:addBS}}
We compare the two unions of cusps from Remark~\ref{rem:addBS}.
As before, we use the procedure \texttt{bernstein} from \texttt{gmssing.lib}, and we recall
that $\lct_0$ is the negative of the largest root of the local Bernstein--Sato polynomial.

\begin{verbatim}
LIB "gmssing.lib";                // local Bernstein--Sato polynomial
ring r = 0,(x,y),ds;              // local ordering at the origin
proc check(poly f)
{
  int mu  = vdim(std(jacob(f)));              // Milnor number
  int tau = vdim(std(jacob(f) + ideal(f)));   // Tjurina number
  list L  = bernstein(f);                     // roots of b_{f,0}(s)
  number lct = -leadcoef(L[1][size(L[1])]);   // largest root is -lct
  printf("mu = %s, tau = %s, lct = %s", mu, tau, lct);
  L;
}
check((y2 - x3)*(y2 - 2x3));      // two tangent cusps
check((x2 - y3)*(y2 - x3));       // two transverse cusps
\end{verbatim}

For the two tangent cusps $(y^2-x^3)(y^2-2x^3)=0$, the output is $\mu_0=\tau_0=15$ and
$\lct_0=5/12$. 

For the two transverse cusps $(x^2-y^3)(y^2-x^3)=0$, the output is $\mu_0=11$, $\tau_0=10$
and $\lct_0=1/2$.


\begin{thebibliography}{99}

\bibitem{Al}
P.~Aluffi, Chern classes of free hypersurface arrangements, \emph{J.~Singul.} \textbf{5} (2012), 19--32.

\bibitem{BHPV}
W.~Barth, K.~Hulek, C.~Peters, and A.~Van de Ven, \emph{Compact complex surfaces}, 2nd ed., Ergeb.\ Math.\ Grenzgeb.\ (3) \textbf{4}, Springer, Berlin, 2004.

\bibitem{BGM}
J.~Brian\c{c}on, M.~Granger, and Ph.~Maisonobe, Sur le polyn\^ome de Bernstein des singularit\'es semi-quasi-homog\`enes, Pr\'epublication de l'Universit\'e de Nice, no.~138, 1986.

\bibitem{BK}
E.~Brieskorn and H.~Kn\"orrer, \emph{Plane Algebraic Curves}, translated from the German by J.~Stillwell,
Modern Birkh\"auser Classics, Birkh\"auser/Springer, Basel, 2012.

\bibitem{BCP}
L.~Bus\'e, F.~Catanese, and E.~Postinghel, \emph{Algebraic curves and surfaces. A history of shapes}, SISSA Springer Ser.\ \textbf{4}, Springer, Cham, 2023.

\bibitem{Sing}
W.~Decker, G.-M.~Greuel, G.~Pfister, and H.~Sch\"onemann, \textsc{Singular} 4.3.2 --- A computer algebra system for polynomial computations, \url{https://www.singular.uni-kl.de}.

\bibitem{Di17}
A.~Dimca, \emph{Hyperplane arrangements. An introduction}, Universitext, Springer, Cham, 2017.

\bibitem{Di24}
A.~Dimca, On free and plus-one generated curves arising from free curves by addition-deletion of a line. \textit{Bull. Math. Sci.} \textbf{14, No. 3} (2024), Article ID 2450007.

\bibitem{D24}
A.~Dimca, On the Castelnuovo--Mumford regularity of curve arrangements,
\textit{Rev. Roum. Math. Pures Appl.} \textbf{70, No. 1-2} (2025), 73--84.

\bibitem{Di25}
A.~Dimca, Some remarks on plane curves related to freeness, \emph{Bull.\ Math.\ Soc.\ Sci.\ Math.\ Roumanie (N.S.)} \textbf{68(116)} (2025), 281--292.

\bibitem{DS}
A.~Dimca and E.~Sernesi, Syzygies and logarithmic vector fields along plane curves, \emph{J.~\'Ec.\ polytech.\ Math.} \textbf{1} (2014), 247--267.

\bibitem{DS17}
A. Dimca and G. Sticlaru, Free divisors and rational cuspidal plane curves, \textit{Math. Res. Lett.} \textbf{24(4)} (2017), 1023--1042.

\bibitem{dPW}
A.\,A.~du Plessis and C.\,T.\,C.~Wall, Application of the theory of the discriminant to highly singular plane curves, \emph{Math.\ Proc.\ Cambridge Philos.\ Soc.} \textbf{126} (1999), 259--266.

\bibitem{ELSV}
L.~Ein, R.~Lazarsfeld, K.\,E.~Smith, and D.~Varolin, Jumping coefficients of multiplier ideals, \emph{Duke Math.\ J.} \textbf{123} (2004), 469--506.

\bibitem{Ka}
M.~Kashiwara, $B$-functions and holonomic systems. Rationality of roots of $B$-functions, \emph{Invent.\ Math.} \textbf{38} (1976/77), 33--53.

\bibitem{Ko}
J.~Koll\'ar, Singularities of pairs, in: \emph{Algebraic geometry---Santa Cruz 1995}, Proc.\ Sympos.\ Pure Math.\ \textbf{62}, Part~1, Amer.\ Math.\ Soc., Providence, RI, 1997, 221--287.

\bibitem{Li}
X.~Liao, Chern classes of logarithmic vector fields, \emph{J.~Singul.} \textbf{5} (2012), 109--114.

\bibitem{Lich}
B.~Lichtin, Poles of $|f(z,w)|^{2s}$ and roots of the $b$-function, \emph{Ark.\ Mat.} \textbf{27} (1989), 283--304.

\bibitem{Ma}
B.~Malgrange, Polyn\^omes de Bernstein--Sato et cohomologie \'evanescente, in: \emph{Analysis and topology on singular spaces, II, III (Luminy, 1981)}, Ast\'erisque \textbf{101--102}, Soc.\ Math.\ France, Paris, 1983, 243--267.

\bibitem{Mi}
J.~Milnor, \emph{Singular points of complex hypersurfaces}, Ann.\ of Math.\ Stud.\ \textbf{61}, Princeton Univ.\ Press, Princeton, NJ, 1968.

\bibitem{MO}
J.~Milnor and P.~Orlik, Isolated singularities defined by weighted homogeneous polynomials, \emph{Topology} \textbf{9} (1970), 385--393.

\bibitem{MP}
M.~Musta\c{t}\u{a} and M.~Popa, Hodge ideals for $\mathbb{Q}$-divisors, $V$-filtration, and minimal exponent, \emph{Forum Math.\ Sigma} \textbf{8} (2020), Paper No.~e19.

\bibitem{Po23}
P. Pokora, $\mathcal{Q}$-conic arrangements in the complex projective plane, \textit{Proc. Am. Math. Soc.} \textbf{151(7)} (2023), 2873--2880.
\bibitem{Po24}
P.~Pokora, On the addition technique for Betti and Poincar\'e polynomials of plane curves, \emph{Bull.\ Pol.\ Acad.\ Sci.\ Math.} \textbf{72} (2024), 111--118.

\bibitem{Po24a}
P. Pokora, On free and nearly free arrangements of conics admitting certain ADE singularities, \textit{Ann. Univ. Ferrara, Sez. VII, Sci. Mat.} \textbf{70(3)} (2024), 593--606.

\bibitem{Po25}
P.~Pokora, On Poincar\'e polynomials for plane curves with quasi-homogeneous singularities, \emph{Bull.\ Lond.\ Math.\ Soc.} \textbf{57} (2025), 2549--2560.

\bibitem{PoBS}
M.~Popa, \emph{The Bernstein--Sato polynomial}, lecture notes. \url{https://people.math.harvard.edu/~mpopa/notes/Bernstein-Sato-notes.pdf}

\bibitem{Sa71}
K.~Saito, Quasihomogene isolierte Singularit\"aten von Hyperfl\"achen, \emph{Invent.\ Math.} \textbf{14} (1971), 123--142.

\bibitem{SaM}
M.~Saito, On microlocal $b$-function, \emph{Bull.\ Soc.\ Math.\ France} \textbf{122} (1994), 163--184.

\bibitem{SaH}
M.~Saito, Bernstein--Sato polynomials of hyperplane arrangements, \emph{Selecta Math.\ (N.S.)} \textbf{22} (2016), 2017--2046.

\bibitem{SaBP}
M.~Saito, Bernstein--Sato polynomials of semi-weighted-homogeneous polynomials of nearly Brieskorn--Pham type, arXiv:2210.01028.

\bibitem{ST}  H.~Schenck, S.~Tohaneanu, Freeness of conic-line arrangements in $\mathbb{P}^2$. 
\textit{Comment. Math. Helv.} {\bf 84(2)} (2009), 235--258.

\bibitem{Wa}
C.\,T.\,C.~Wall, \emph{Singular points of plane curves}, London Math.\ Soc.\ Stud.\ Texts \textbf{63}, Cambridge Univ.\ Press, Cambridge, 2004.

\bibitem{Wal}
U.~Walther, Bernstein--Sato polynomial versus cohomology of the Milnor fiber for generic hyperplane arrangements, \emph{Compos.\ Math.} \textbf{141} (2005), 121--145.
\end{thebibliography}
\end{document}